\documentclass[a4paper,9pt]{article}
\usepackage{amsfonts}
\usepackage{amssymb,amsthm,color,tikz}
\usepackage{amsmath}
\usepackage{comment}

\usepackage[normalem]{ulem}

\newcommand{\R}{\mathbb{R}}

\newtheorem{theorem}{Theorem}
\newtheorem{lemma}[theorem]{Lemma}
\newtheorem{proposition}[theorem]{Proposition}

\def\cp{\mathrm{cap}\,}

\newcommand{\bp}{\begin{proof}}
\newcommand{\ep}{\end{proof}}

\begin{document}
\title{On the torsion function with mixed boundary conditions}

\author{{M. van den Berg} \\
School of Mathematics, University of Bristol\\
Fry Building, Woodland Road, Bristol BS8 1UG\\
United Kingdom\\
\texttt{mamvdb@bristol.ac.uk}\\
\\
Tom Carroll\\
Department of Mathematics\\
University College Cork\\
Cork, Ireland\\
\texttt{t.carroll@ucc.ie}}
\date{14 June 2020}\maketitle
\vskip 3truecm \indent
\begin{abstract}\noindent
Let $D$ be a non-empty open subset of $\R^m,\,m\ge 2$, with boundary $\partial D$, with finite Lebesgue measure $|D|$, and which satisfies a parabolic Harnack principle. Let $K$ be a compact, non-polar subset of $D$. We obtain the leading asymptotic behaviour as $\varepsilon\downarrow 0$ of the $L^{\infty}$ norm of the torsion function with a Neumann boundary condition on
$\partial D$, and a Dirichlet boundary condition on $\partial (\varepsilon K)$, in terms of the first eigenvalue of the Laplacian with corresponding boundary conditions.
These estimates quantify those of Burdzy, Chen and Marshall who showed that $D\setminus K$ is a non-trap domain.
\end{abstract}
\begin{center} \textbf{Keywords}: Torsion function, Dirichlet boundary condition, Neumann boundary condition.
\end{center}


\section{Introduction and main results \label{sec1}}
Let $D$ be an open, non-empty set in $\R^m,\,m\ge2$, with finite Lebesgue measure $|D|$, and let $K\subset D$ be a compact set with boundary $\partial K$, and with positive logarithmic capacity if $m=2$ or with positive Newtonian capacity $\cp(K)$ if $m\ge 3$.
Let $u_{K,D}$ be the solution of
\begin{equation*}
-\Delta u=1,
\end{equation*}
with Dirichlet boundary condition
\begin{equation}\label{e2}
u(x)=0,\, x\in \partial K,
\end{equation}
and Neumann boundary condition
 \begin{equation}\label{e3}
\frac{\partial u}{\partial \nu}(x)=0,\,x\in\partial D,
\end{equation}
where $\nu$ is the inward normal.
Boundary conditions \eqref{e2} and \eqref{e3} have to be understood in the weak sense.
In particular \eqref{e2} holds for all regular points of $\partial K$.
Let $\pi_D(x,y;t)$, $x\in D$, $y\in D$, $t>0$ denote the Neumann heat kernel for $D$. We say that the parabolic Harnack principle (PHP for short) holds in $D$ if for some $t_0\in(0,\infty)$ there exists $c_0=c_0(D,t_0)<\infty$, such that
\begin{equation*}
\pi_D(x,y;t)\le c_0\pi_D(v,w;t),\, t\ge t_0,\,\,x,y,v,w\in D.
\end{equation*}
See also \cite{LSC}.
As was pointed out in \cite{BCM}, PHP is equivalent to the following assertion: there exist $t_1\in (0,\infty)$, $c_1<\infty$, $c_2>0$ depending on $D$ such that
\begin{equation}\label{e5}
\sup_{x,y\in D}\bigg|\pi_D(x,y;t)-\frac{1}{|D|}\bigg|\le c_1e^{-c_2t},\,t\ge t_1.
\end{equation}

It was shown in \cite{BCM} that if $D$ satisfies PHP then $u_{K,D}$ is bounded, and $D\setminus K$ is a {\it non-trap} domain.
In Theorem \ref{the1} below we quantify this statement in terms of the first eigenvalue $\lambda(K,D)$ of the Laplacian with boundary conditions \eqref{e2} and \eqref{e3} in the case where $K$ is scaled down by a factor $\varepsilon$ with respect to
a fixed point (the origin) in $D$.

Estimates of this type are well known for the torsion function $u_{\Omega}$ for an open set $\Omega$ satisfying a $0$ Dirichlet boundary condition on $\partial \Omega$. In
\cite{vdBC} it was shown that $u_{\Omega}\in L^{\infty}(\Omega)$ if and only if $\lambda(\Omega)>0$. If the latter holds then
\begin{equation*}
\lambda(\Omega)^{-1}\le \|u_{\Omega}\|_{\infty}\le \mathfrak{c}_m\lambda(\Omega)^{-1},
\end{equation*}
where $\mathfrak{c}_m$ is the sharp constant defined by
\begin{equation*}
\mathfrak{c}_m=\sup\{\lambda(\Omega)\|u_{\Omega}\|_{\infty}:\, \Omega\, \textup{open in}\, \R^m, \lambda(\Omega)>0\},
\end{equation*}
and $\|\cdot\|_p$ denotes the standard $L^p$ norm, $1\le p\le \infty$.

In \cite{vdBC} it was shown that $\mathfrak{c}_m\le4+3m\log 2$. This bound has been improved since. See for example \cite{GS} and \cite{HV}.
For general open, non-empty, and connected $D$, and a non-empty compact subset $K\subset D$ one does not have boundedness of $u_{K,D}$. Examples of these {\it trap} domains were given in \cite{BCM}.
\begin{theorem}\label{the1}
Let $D\subset \R^m$, $m\ge 2$, be open, non-empty, containing the origin, and let $D$ satisfy the parabolic Harnack principle. If  $K$ is a non-polar compact subset of $D$, then for $\varepsilon\downarrow 0$,
\begin{equation}\label{e8}
\lambda(\varepsilon K,D)\|u_{\varepsilon K,D}\|_{\infty}=\begin{cases}1+O\big((\log\varepsilon^{-1})^{-1/2}\big),\,\,m=2,\\ 1+O\big(\varepsilon^{(m-2)/2}\big),\,\,\,m\ge 3,  \end{cases}
\end{equation}
where $\varepsilon K=\{y\in \R^m:\varepsilon^{-1}y\in K\}$.
Furthermore for any non-polar compact set $K\subset D$,
\begin{equation}\label{e9}
\|u_{ K,D}\|_{\infty}\ge \frac{1}{\lambda(K,D)}.
\end{equation}
\end{theorem}

It was shown in Theorem 2.5(i) in \cite{BCM} that if \eqref{e5} holds, then the Neumann Laplacian on $D$ has discrete spectrum. Sufficient geometric conditions for $D$ to satisfy the PHP were obtained in, for example, Corollary 2.7 of \cite{BCM}. Conversely PHP implies some geometric and spectral properties of $D$. The proposition below is of independent interest.
\begin{proposition}\label{prop1}
Let $D$ be open, non-empty, with $|D|<\infty$. If \eqref{e5} holds then we have the following.
\begin{enumerate}
\item[\textup{(i)}] $D$ is connected.
\item[\textup{(ii)}] The first eigenvalue of the Neumann Laplacian acting in $L^2(D)$ has multiplicity $1$.
\item[\textup{(iii)}]
\begin{equation}\label{e10}
\mu(B)\bigg(\frac{|B|}{|D|}\bigg)^{2/m}\ge \mu(D)\ge c_2,
\end{equation}
where $\mu(D)$ is the first non-zero eigenvalue of the Neumann Laplacian acting in $L^2(D)$, and $B$ is a ball of radius $1$ in $\R^m$.
\end{enumerate}
\end{proposition}

\section{Proof of Theorem \ref{the1} \label{sec2}}
In this section we prove Theorem \ref{the1}.
\begin{proof}
Let $\pi_{K,D}(x,y;t),\,x\in D\setminus K$, $y\in D\setminus K$, $t>0$ denote the heat kernel with a Neumann boundary condition on $\partial D$, and with a $0$ Dirichlet boundary condition on $\partial K$.
We have  for $\delta\in (0,1)$,
\begin{align}\label{e11}
u_{K,D}(x)&=\int_0^{\infty}dt\,\int_{D\setminus K}dy\,\pi_{K,D}(x,y;t)\nonumber \\ &
=\int_0^{t_1/(1-\delta)}dt\,\int_{D\setminus K}dy\,\pi_{K,D}(x,y;t)+\int_{t_1/(1-\delta)}^{\infty}dt\,\int_{D\setminus K}dy\,\pi_{K,D}(x,y;t)\nonumber \\ &
\le \int_0^{t_1/(1-\delta)}dt\,\int_{D\setminus K}dy\,\pi_{D}(x,y;t)+\int_{t_1/(1-\delta)}^{\infty}dt\,\int_{D\setminus K}dy\,\pi_{K,D}(x,y;t)\nonumber \\ &
\le \frac{t_1}{1-\delta}+\int_{t_1/(1-\delta)}^{\infty}dt\,\int_{D\setminus K}dy\,\pi_{K,D}(x,y;t).
\end{align}
By the heat semigroup property, and by Cauchy-Schwarz's inequality,
\begin{align}\label{e12}
\pi_{K,D}(x,y;t) & =   \int_{D\setminus K} \pi_{K,D}(x,z;t/2)\, \pi_{K,D}(z,y;t/2)\,dz \cr & \leq
\left( \int_{D\setminus K} \pi_{K,D}(x,z;t/2)^2\,dz\right)^{1/2} \left( \int_{D\setminus K}
\pi_{K,D}(z,y;t/2)^2\,dz\right)^{1/2}\cr & =   \big( \pi_{K,D}(x,x;t)\,
\pi_{K,D}(y,y;t) \big)^{1/2}.
\end{align}
By the spectral theorem we have
\begin{equation}\label{e13}
\pi_{K,D}(x,x;t)\le e^{-\delta t\lambda(K,D)}\pi_{K,D}(x,x;(1-\delta)t).
\end{equation}
By \eqref{e12} and \eqref{e13},
\begin{align}\label{e14}
\big(\pi_{K,D}(x,y;t)\big)^{\delta}&\le e^{-\delta^2t\lambda(K,D)}\big(\pi_{K,D}(x,x;(1-\delta)t)\pi_{K,D}(y,y;(1-\delta)t)\big)^{\delta/2}\nonumber \\ &\le
e^{-\delta^2 t\lambda(K,D)}\sup_{x,y\in D}\big(\pi_{K,D}(x,y;(1-\delta)t)\big)^{\delta}\nonumber \\ &\le
e^{-\delta^2t\lambda(K,D)}\sup_{x,y\in D}\big(\pi_{D}(x,y;(1-\delta)t)\big)^{\delta}.
\end{align}
By \eqref{e5},
\begin{align*}
\big(\pi_{D}(x,y;(1-\delta)t)\big)^{\delta}&\le \bigg(\frac{1}{|D|}+c_1e^{-c_2(1-\delta)t}\bigg)^{\delta}\nonumber \\ &
\le \frac{1}{|D|^{\delta}}+c_1^{\delta}e^{-c_2\delta(1-\delta) t},\,t\ge \frac{t_1}{1-\delta}.
\end{align*}
This, together with \eqref{e14}, gives
\begin{equation}\label{e16}
\big(\pi_{K,D}(x,y;t)\big)^{\delta}\le e^{-\delta^2t\lambda(K,D)}\bigg(\frac{1}{|D|^{\delta}}+c_1^{\delta}e^{-c_2\delta(1-\delta)t}\bigg),\,t\ge \frac{t_1}{1-\delta}.
\end{equation}
We obtain by \eqref{e16}, and  by H\"older's inequality,
\begin{align}\label{e17}
&\int_{t_1/(1-\delta)}^{\infty}dt\,\int_{D\setminus K}dy\,\pi_{K,D}(x,y;t)\nonumber \\ &\hspace{2mm}
\le \int_{t_1/(1-\delta)}^{\infty}dt\,\int_{D\setminus K}dy\,\big(\pi_{K,D}(x,y;t)\big)^{1-\delta}e^{-\delta^2t\lambda(K,D)}\bigg(\frac{1}{|D|^{\delta}}+c_1^{\delta}e^{-c_2\delta(1-\delta)t}\bigg)\nonumber \\ &\hspace{2mm}
\le \int_{t_1/(1-\delta)}^{\infty}dt\,\int_{D}dy\,\big(\pi_{D}(x,y;t)\big)^{1-\delta}e^{-\delta^2t\lambda(K,D)}\bigg(\frac{1}{|D|^{\delta}}+c_1^{\delta}e^{-c_2\delta(1-\delta)t}\bigg)\nonumber \\ &\hspace{2mm}
\le \int_{t_1/(1-\delta)}^{\infty}dt\,\bigg(\int_{D}dy\,\pi_{D}(x,y;t)\bigg)^{1-\delta}|D|^{\delta}e^{-\delta^2 t \lambda(K,D)}\bigg(\frac{1}{|D|^{\delta}}+c_1^{\delta}e^{-c_2\delta(1-\delta)t}\bigg)\nonumber \\ &\hspace{2mm}
=\frac{1}{\delta^2\lambda(K,D)}e^{-\delta^2{t_1}\lambda(K,D)/(1-\delta)}\nonumber \\ & \hspace{23mm}+c_1^{\delta}|D|^{\delta}\big(c_2\delta(1-\delta)+\delta^2\lambda(K,D)\big)^{-1}e^{-{t_1}(\delta c_2+\delta^2\lambda(K,D)/(1-\delta))}\nonumber \\&\hspace{2mm}
\le \frac{1}{\delta^2\lambda(K,D)}+\frac{c_1^{\delta}|D|^{\delta}}{c_2\delta(1-\delta)}.
\end{align}
By \eqref{e11} and \eqref{e17},
\begin{equation*}
u_{K,D}(x)\lambda(K,D)\le \delta^{-2}+\bigg(\frac{t_1}{1-\delta}+\frac{c_1^{\delta}|D|^{\delta}}{c_2\delta(1-\delta)}\bigg)\lambda(K,D).
\end{equation*}
By taking the supremum over all $x\in D\setminus K$ we obtain
\begin{equation*}
\|u_{K,D}\|_{\infty}\lambda( K,D)\le \delta^{-2}+\bigg(\frac{t_1}{1-\delta}+\frac{c_1^{\delta}|D|^{\delta}}{c_2\delta(1-\delta)}\bigg)\lambda(K,D).
\end{equation*}
Hence for $\delta\in(0,1)$ and $\varepsilon\in(0,1)$,
\begin{align}\label{e20}
\|u_{\varepsilon K,D}\|_{\infty}\lambda(\varepsilon K,D)&\le \delta^{-2}+\bigg(\frac{t_1}{1-\delta}+\frac{c_1^{\delta}|D|^{\delta}}{c_2\delta(1-\delta)}\bigg)\lambda(\varepsilon K,D).
\end{align}
In the lemma below we obtain an upper bound for the rate at which $\lambda(\varepsilon K,D)\downarrow 0$ as $\varepsilon\downarrow 0$.

\begin{lemma}\label{lem1}
If $D$ is open, non-empty in $\R^m,\,m\ge 3$, with $|D|<\infty$, and if $K\subset D$ with $\cp(K)>0$
then
\begin{equation}\label{e21}
\limsup_{\varepsilon\downarrow 0}\varepsilon^{2-m}\lambda(\varepsilon K,D)\le \frac{\cp(K)}{|D|}.
\end{equation}
If $D$ is open, non-empty in $\R^2$, with $|D|<\infty$, and if $K\subset D$ has strictly positive logarithmic capacity, then
\begin{equation}\label{e22}
\limsup_{\varepsilon\downarrow 0}\big(\log \varepsilon^{-1}\big)\lambda(\varepsilon K,D)\le \frac{2\pi}{|D|}.
\end{equation}
\end{lemma}
We note that (i) the constants in the right-hand sides of \eqref{e21} and \eqref{e22} are well-known and sharp (see for example \cite{SO}), (ii) both formulae hold for arbitrary open and connected sets $D$ with $|D|<\infty$, and without any regularity assumptions on $\partial D$.
We now choose
\begin{equation}\label{e23}
\delta=1-|D|^{1/m}\lambda(\varepsilon K,D)^{1/2}.
\end{equation}
Then $\delta\in(0,1)$ for all $\varepsilon$ sufficiently small. By \eqref{e20} and \eqref{e23},
\begin{equation}\label{e24}
\|u_{\varepsilon K,D}\|_{\infty}\lambda(\varepsilon K,D)\le1+O\big(\lambda(\varepsilon K,D)^{1/2}\big).
\end{equation}

The proof of \eqref{e9} is similar to the one of Theorem 5 in \cite{vdB}, and Theorem 1, (0.5) in \cite{BC}. Let $\psi$ denote the normalised first eigenfunction (positive) of the Laplacian with Neumann and Dirichlet boundary conditions on $\partial D$ and $\partial K$ respectively, suppressing both $K$ and $D$ dependence.
We have by Cauchy-Schwarz's inequality that $\int_{D\setminus K}\psi\le |D\setminus K|^{1/2}.$ Using
\begin{equation*}
\psi \frac{\partial u_{K,D}}{\partial \nu}=u_{K,D}\frac{\partial \psi}{\partial \nu}=0 \,\, \textup{on}\, \partial D\cup \partial K,
\end{equation*}
we obtain by Green's formula,
\begin{align*}
\lambda(K,D)\|u_{K,D}\|_{\infty}\int_{D\setminus K}\psi&\ge \lambda(K,D)\int_{D\setminus K}u_{K,D}\psi =-\int_{D\setminus K}u_{K,D}\Delta\psi\nonumber \\ &=-\int_{D\setminus K}\psi \Delta u_{K,D}= \int_{D\setminus K}\psi.
\end{align*}
This implies the assertion.

Finally \eqref{e8} follows by \eqref{e9}, \eqref{e24}, and Lemma \ref{lem1}.
\end{proof}

\section{Proof of Lemma \ref{lem1} and Proposition \ref{prop1} \label{sec3}}
\noindent{\it Proof of Lemma \textup{\ref{lem1}}.}
Recall that $0\in D$, and so
\begin{equation*}
R= \min\{|y|:y\in \partial D\}>0.
\end{equation*}
Since $K$ is compact,
\begin{equation*}
R_K=\max\{|x|:x\in K\}<\infty.
\end{equation*}
Let
\begin{equation*}
\varepsilon_1=\min\bigg\{1,\frac{R}{R_K}\bigg\}.
\end{equation*}
If $\varepsilon\le \varepsilon_1$ then $\varepsilon K\subset B(0;R)$.
See \cite{MET} for estimates related to the proof of Lemma \ref{lem1}.
First we consider the case $m\ge 3$. Let $\mu_K$ denote the equilibrium measure of $K$ in $\R^m$, and let
\begin{equation*}
\phi_K(x)=\frac{\Gamma((m-2)/2)}{4\pi^{m/2}}\int_K\mu_K(dy)\,|x-y|^{2-m}.
\end{equation*}
Then $\phi_K(x)=1,\,x\in K$, $0<\phi_K<1,\,x\in \R^m\setminus K$, and $\phi_K$ is smooth on the complement of $K$. We use $1-\phi_K$ as a trial function in the Rayleigh-Ritz characterisation
of $\lambda(K,D)$. This gives
\begin{align}\label{e31}
\lambda(K,D)&=\inf_{u\in H^1(D),\,u|_{K}=0}\frac{\int_{D\setminus K}|\nabla u|^2}{\int_{D\setminus K}u^2}\nonumber \\ &
\le \frac{\int_{D\setminus K}|\nabla \phi_K|^2}{\int_{D\setminus K}(1-\phi_K)^2}\nonumber \\ &\le
\frac{\int_{\R^m\setminus K}|\nabla\phi_K|^2}{\int_{D\setminus K}(1-\phi_K)^2}\nonumber \\ &
=\frac{\cp(K)}{\int_{D\setminus K}(1-\phi_K)^2}.
\end{align}
It remains to bound the denominator in the right-hand side of \eqref{e31} from below.
Since we will apply this lower bound with $\varepsilon_1K$ rather than $K$ itself,
we assume that $K \subset B(0;R)$. We let $0<\alpha<1$.
It is a standard fact that the capacitary potential is monotone increasing in $K$. In particular,
\begin{equation*}
\phi_K(x)\le \phi_{B(0;R)}(x)=\min\bigg\{1,\bigg(\frac{R}{|x|}\bigg)^{m-2}\bigg\}.
\end{equation*}
Hence

\begin{align}\label{e33}
\int_{D\setminus K}(1-\phi_K)^2&
\ge(1-\alpha)^2\int_{\{\phi_K(x)\le\alpha\}\cap D} 1\nonumber \\ &
\ge (1-\alpha)^2 \big(|D|-|\{\phi_{B(0;R)}(x)>\alpha\}|\big)\nonumber \\ &
\ge (1-\alpha)^2 \big(|D|-\alpha^{-m/(m-2)}\omega_mR^m\big),
\end{align}
where $\omega_m=|B_1(0)|$. We choose $\alpha$ such that
\begin{equation}\label{e34}
\alpha=\alpha^{-m/(m-2)}\frac{|B(0;R)|}{|D|}.
\end{equation}
This, together with \eqref{e31}, \eqref{e33}, and \eqref{e34} implies
\begin{equation}\label{e35}
\lambda(K,D)\le \frac{\cp(K)}{|D|}\bigg(1-\bigg(\frac{|B(0;R)|}{|D|}\bigg)^{(m-2)/(2(m-1))}\bigg)^{-3}.
\end{equation}
In particular for $\varepsilon\in (0,1]$, $\varepsilon \varepsilon_1 K\subseteq \varepsilon B(0;R)$,
and this together with \eqref{e35} gives
\begin{equation}\label{e36}
\lambda(\varepsilon \varepsilon_1 K,D)\le \frac{\cp(\varepsilon \varepsilon_1K)}{|D|}\bigg(1-\bigg(\frac{\varepsilon|B(0;R)|}{|D|}\bigg)^{(m-2)/(2(m-1))}\bigg)^{-3}.
\end{equation}
Formula \eqref{e21} follows by \eqref{e36}, and scaling of the Newtonian capacity,
\begin{equation*}
\cp(\varepsilon K)=\varepsilon^{m-2} \cp(K).
\end{equation*}

Next we consider the planar case $m=2$. We use Hadamard's method of descent so as to avoid logarithmic potential theory. See for example p.51 in \cite{MET}.
Let $h\ge R$, and consider the cylinder $(D\setminus K)\times (0,h)\subset \R^3$. Then the first eigenvalue of the Laplacian acting in $L^2(D\setminus K))$ with
Dirichlet boundary condition on $\partial K$, and Neumann boundary condition on $\partial D$ is precisely equal to
the first eigenvalue of the Laplacian acting in $L^2((D\setminus K)\times (0,h))$ with Dirichlet boundary condition on $\partial (K\times(0,h))$, and Neumann boundary condition on $\partial (D\times (0,h))\setminus \partial (K\times(0,h)).$
We apply \eqref{e35} to the setting above and obtain by monotonicity of Newtonian capacity,
\begin{align}\label{e38}
\lambda(\varepsilon \varepsilon_1K,D)&\le \lambda(\varepsilon B(0;R),D)\nonumber\\&
\le \frac{\cp( B(0;\varepsilon R)\times(0,h))}{|D|h}\bigg(1-\bigg(\frac{\varepsilon|B(0;R)|}{|D|}\bigg)^{1/4}\bigg)^{-3}.
\end{align}
To obtain an upper bound on $\cp( B(0;\varepsilon R)\times(0,h))$ we let $C(R',h')\subset \R^3$ be an ellipsoid with a circular cross section of radius $R'$ and axis $h'$. Then for a suitable translation and rotation
$C(R',h')\supset B(0;\varepsilon R)\times(0,h)$ provided
\begin{equation}\label{e39}
\frac{h^2}{h'^2}+\frac{(\varepsilon R)^2}{R'^2}\le 1.
\end{equation}
We let $\alpha\in (0,1)$ be arbitrary, and choose
\begin{equation}\label{e40}
R'=\varepsilon^{-\alpha}(\varepsilon R),
\end{equation}
and
\begin{equation}\label{e41}
h'=\big(1-\varepsilon^{2\alpha}\big)^{-1/2}h.
\end{equation}
The choice \eqref{e40}--\eqref{e41} satisfies \eqref{e39}. For $\frac{h'}{R'}\rightarrow\infty$, or equivalently $\varepsilon\downarrow 0$ with $h$ fixed, we have by formula (12) on p.260 in \cite{IMcK},
\begin{align*}
\cp(C(R',h'))&=\frac{2\pi h'}{\log(h'/R')}(1+o(1))\nonumber \\ &
\le \frac{2\pi h}{\big(1-\varepsilon^{2\alpha}\big)^{1/2}\log(h/R')}(1+o(1))\nonumber \\ &
\le \frac{2\pi h}{(1-\alpha)\big(1-\varepsilon^{2\alpha}\big)^{1/2}\log\varepsilon^{-1}}(1+o(1)).
\end{align*}
Thus,
\[
\frac{\cp( B(0;\varepsilon R)\times(0,h))}{|D|h} \leq \frac{2\pi }{(1-\alpha) \vert D \vert
	\log\varepsilon^{-1}}(1+o(1)).
\]
By \eqref{e38},
\[
\limsup_{\varepsilon\downarrow 0}\big(\log \varepsilon^{-1}\big)\lambda(\varepsilon \varepsilon_1 K,D)
\leq \frac{2\pi }{(1-\alpha)\,\vert D \vert}.
\]
Since $\alpha\in (0,1)$ was arbitrary, this completes the proof of the case $m=2$.

{\hspace*{\fill }$\square $}
\medskip

\noindent {\it Proof of Proposition \textup{\ref{prop1}}.}
To prove (i) we recall that, since $D$ is open, $D$ is a countable union of open components. Suppose that this union contains at least two elements, one of which is $C$. Then both $C$ and $D\setminus C$ are open and non-empty.
Let $1_A$ denote the indicator function of a set $A$. From \eqref{e5} we obtain,
\begin{equation*}
\bigg|\int _Cdy\, \pi_D(x,y;t)-\frac{|C|}{|D|}\bigg|\le c_1|C|e^{-c_2t},\, t\ge t_1,\, x\in D.
\end{equation*}
We note that
\begin{equation*}
q_{C,D}(x;t)=\int _Cdy\, \pi_D(x,y;t)
\end{equation*}
is the solution of the heat equation
\begin{equation*}
\Delta q=\frac{\partial q}{\partial t},
\end{equation*}
with initial condition
\begin{equation*}
q(x;0)=1_C(x),
\end{equation*}
and with a Neumann (insulating) boundary condition on $\partial D$. It follows that
\begin{equation*}
q_{C,D}(x;t)=1_C(x),\,t>0.
\end{equation*}
From \eqref{e5} we have
\begin{equation*}
\bigg|1-\frac{|C|}{|D|}\bigg|\le c_1|C|e^{-c_2t},\, t\ge t_1,\, x\in C.
\end{equation*}
We conclude that, by taking the limit $t\rightarrow \infty$, $|C|=|D|$. Since $C\subset D$, $|D\setminus C|=0$. This contradicts $D\setminus C$ is open and non-empty. This in turn implies that $D$ consists of just one component
$C$. Hence $C$ is connected. This implies assertion (ii).
To prove (iii) we have that \eqref{e5} implies
\begin{equation*}
\int_Ddx\,\pi_D(x,x;t)\le 1+ c_1|D|e^{-c_2t},\, t\ge t_1.
\end{equation*}
Hence the Neumann heat semigroup is trace-class, and
\begin{equation}\label{e52}
1+e^{-t\mu(D)}\le \int_Ddx\,\pi_D(x,x;t)\le 1+ c_1|D|e^{-c_2t},\, t\ge t_1.
\end{equation}
Taking the limit $t\rightarrow\infty$ in \eqref{e52} implies the second inequality in \eqref{e10}. The first inequality in \eqref{e10} is due to Weinberger \cite{W}.
{\hspace*{\fill }$\square $}\\

\medskip
\section*{Acknowledgements} MvdB acknowledges support by the Leverhulme Trust through Emeritus Fellowship EM-2018-011-9.

\end{document}